\documentclass[pdftex, a4paper, 12pt, oneside]{amsart}
\usepackage{geometry}
\usepackage[cp1250]{inputenc}

\usepackage{graphicx}
\usepackage{lpic}
\usepackage{url}
\usepackage{amsfonts, amstext, amsmath, amsthm, amscd, amssymb}
\usepackage[sc]{mathpazo}
\usepackage{enumitem, textcomp, setspace}
\usepackage{caption}

\usepackage{tikz}
\usepackage{tikz-cd}
\usepackage{nicefrac}
\usepackage{longtable}

\numberwithin{equation}{section}

\theoremstyle{plain}
\newtheorem{theorem}{Theorem}[section]
\newtheorem{proposition}[theorem]{Proposition}

\newtheorem{conjecture}[theorem]{Conjecture}

\theoremstyle{definition}

\newtheorem{question}[theorem]{Question}

\definecolor{newblue}{rgb}{0.27, 0.32, 0.86}
\definecolor{newred}{rgb}{0.86, 0.32, 0.27}

\providecommand{\keywords}[1]{\textbf{\textit{Key words and phrases:}} #1}
\providecommand{\subjclass}[1]{\textbf{\textit{2020 Mathematics Subject Classification:}} #1}

\graphicspath{ {} } 

\newcommand{\pic}[2]{\raisebox{-.4\height}{\includegraphics[scale=#2]{#1}}}
\newcommand\Xor{\pic{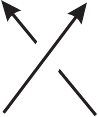}{.50}}
\newcommand\Yor{\pic{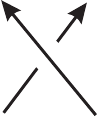} {.50}}
\newcommand\Ior{\pic{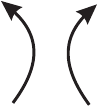} {.50}}

\usepackage[pagebackref=false]{hyperref}

\usepackage[all]{hypcap}

\definecolor{internalLink}{rgb}{0,0,0.5}
\definecolor{citeLink}{rgb}{0,0.5,0}
\definecolor{urlLink}{rgb}{0,0.5,0.5}

\hypersetup{
	linkcolor=internalLink,
	citecolor=citeLink,
	urlcolor=urlLink
}

\begin{document}
	

\title[Inequalities between knot invariants, skein tree depth and delta-crossing numbers]{Integer inequalities between knot invariants,\\ skein tree depth and delta-crossing numbers}

\author{Micha{\l} Jab{\l}onowski}

\address{Institute of Mathematics, Faculty of Mathematics, Physics and Informatics,\newline University of Gda\'nsk, 80-308 Gda\'nsk, Poland}



\email{michal.jablonowski@gmail.com}

\date{\today}

\maketitle


	\begin{abstract}
		The maximum length of the shortest path from a leaf to the root of a skein tree for knots and links gives a measure of the complexity of computing link polynomials by the skein relation (the Jones polynomial, the Alexander-Conway polynomial, and more generally the HOMFLY-PT polynomial). We combine theoretical and computational results on the skein tree depth of knots and links.
		We prove the new upper bound on the skein tree depth of a link and give examples of links where the new bound is stronger than the known bound. We also give the new lower bound. Moreover, we derive tables of knots and links with their skein tree depth that were up to now undetermined (for some of them, we give their range of possible values).
		The paper surveys known (and new) inequalities between integer-valued classical knot invariants. It features a visual graph of the relations.
	\end{abstract}
	
	\section{Introduction}
	A fundamental approach to computing link invariants is via skein relations. Iteratively resolving crossings of a link diagram using a skein relation produces a binary \emph{skein tree} whose root is the original link and whose leaves are unlink diagrams. The \emph{skein tree depth} of a knot or link $L$, denoted $td(L)$, is the minimum depth (length of the longest root-to-leaf path) among all skein resolving trees for $L$. This quantity provides a measure of the complexity of computing link polynomials (such as the Jones, Alexander-Conway, or HOMFLY-PT polynomials) by skein methods. Separately, any knot or link can be described by a \emph{delta-crossing diagram} -- a diagram composed of certain three-strand tangle blocks called \emph{delta-crossings} \cite{NSS15}. The \emph{delta-crossing number} $c_\Delta(K)$ of a knot $K$ is the least number of delta-crossings required in such a diagram \cite{NSS15}. 
	
	In this paper, we explore new relationships between $td(L)$ and other knot invariants, with a focus on $c_\Delta$. We prove a new inequality $2\,c_\Delta(K) \ge td(K)$ for any knot $K$, indicating that resolving a diagram with $n$ delta-crossings requires at least $2n$ skein moves in the worst case. We also present an improved upper bound on $td(L)$ in terms of the braid length and strand numbers in a braid representation of $L$. In addition, we review known bounds relating $td(L)$ to classical invariants like the crossing number and the knot genus. On the computational side, we extend previous calculations of skein tree depth. In particular, building on prior work that classified all links with $td\le 2$ \cite{Tho89}, we determine the exact skein tree depth for all prime knots and links through 8 crossings, and for many links with 9 or more crossings. Our tables include several cases where only a range for $td$ was previously known. We highlight examples where our new inequalities attain equality, and observe that all data so far are consistent with a conjectured relationship between $td(L)$ and the triple-crossing number $c_3(L)$. 
	
	The paper is organized as follows. In Section\;\ref{s2} we give precise definitions of skein tree depth, delta-crossing number, and other knot invariants used later. In Section\;\ref{s3} we present the graphs of selected invariant inequalities. Section\;\ref{s4} contains inequalities involving $td(L)$, including bounds in terms of the crossing number, the oriented genus, and the HOMFLY-PT polynomial. Section\;\ref{s5} introduces and proves the new inequality $2\,c_\Delta(K)\ge td(K)$. Section\;\ref{s6} presents the new upper bound on $td(L)$ derived from braid representations. In Section\;\ref{s7} we discuss computational results and provide tables of $td$ values for knots and links, along with illustrative examples. Finally, Section\;\ref{s8} collects open problems and conjectures related to skein depth and delta-crossing number.
	
	\section{Definitions}\label{s2}
	\subsection{Skein tree depth}
	Any oriented knot or link $L$ in $S^3$ can be reduced to a collection of unlinks by a sequence of skein moves. In particular, given an oriented link diagram $L$, one can perform a \textit{skein resolution} at a crossing: replacing a positive crossing $L_{+}$ by the two diagrams $L_{-}$ (the negative crossing resolution) and $L_{0}$ (the zero or smoothing resolution), or symmetrically replacing a negative crossing $L_{-}$ by $L_{+}$ and $L_{0}$. 
	
	$$ L_+ =\Xor, \quad L_-=\Yor, \quad L_0=\Ior.$$
	
	Iterating this procedure produces a rooted binary tree called a \textit{skein tree}, whose root is $L$ and whose leaves are unlink diagrams. We conventionally draw the $L_{\pm}$ child to the left and the $L_{0}$ child to the right at each branching. The \emph{skein tree depth} of an oriented link $L$, denoted $td(L)$, is defined as the minimum possible depth of a skein tree for $L$. Equivalently, $td(L)$ is the length of the shortest maximal root-to-leaf path in an optimally chosen skein tree for $L$. Intuitively, $td(L)$ measures the minimum number of sequential skein moves required to fully resolve $L$ into unlinks. We call a skein tree achieving this minimum depth a \emph{minimal skein tree} for $L$.
	
	\begin{figure}[h!t]
		\begin{center}
			\begin{lpic}[]{1_tee(9.3cm)}
				\lbl[t]{130,260; $K9n4$}
				\lbl[t]{82,120; $K3a1\#L2a1\; = $}
				\lbl[t]{63,72; $K3a1\#L2a1$}
				\lbl[t]{164,116; $T_1$}
				\lbl[r]{30,45; $K3a1\sqcup T_1$}
				\lbl[l]{159,45; $L2a1$}
				\lbl[l]{150,72; $K3a1$}
				\lbl[l]{52,100; $K6a2$}
				\lbl[b]{40,200; $K8a10$}
				\lbl[b]{145,200; $L7n2$}
				\lbl[r]{33,3; $T_3$}
				\lbl[t]{144,20; $T_2$}
			\end{lpic}
			\caption{A skein tree of the knot $9n4$ with depth $6$.\label{1_tree}}
		\end{center}
	\end{figure}

	For example, Figure\;\ref{1_tree} shows a skein tree for the knot $9_{n4}$ of depth 6. Figure\;\ref{2_tree} shows an alternative skein tree for the same knot (using the same minimal crossing diagram as root) with depth 5. This demonstrates that even for a fixed prime knot and a fixed minimal diagram, different choices of resolution order can yield different depths. It remains an open question whether the skein tree depth of a knot is always realized on a minimal-crossing diagram (see Question~\ref{ques:minimal-diagram} in Section\;\ref{s8}).
	
	\begin{figure}[h!t]
		\begin{center}
			\begin{lpic}[]{2_tee(9.6cm)}
				\lbl[t]{115,220; $K9n4$}
				\lbl[l]{85,80; $K3a1\#L2a1$}
				\lbl[t]{115,130; $T_1$}
				\lbl[r]{28,45; $K3a1\sqcup T_1$}
				\lbl[l]{110,150; $L2a1$}
				\lbl[l]{85,45; $K3a1$}
				\lbl[l]{52,135; $K6a2$}
				\lbl[r]{31,3; $T_3$}
				\lbl[r]{20,20; $T_2$}
			\end{lpic}
			\caption{A skein tree of the knot $9n4$ with depth $5$.\label{2_tree}}
		\end{center}
	\end{figure}

	\subsection{Delta-crossing number}
	A \emph{delta-crossing} is a tangle consisting of three crossing strands (often depicted as a triangular arrangement of crossings). It is known that any knot or link has a diagram that can be decomposed entirely into delta-crossing tangles \cite{NSS15}. The \emph{delta-crossing number} $c_{\Delta}(K)$ of a knot $K$ is defined as the minimum number of delta-crossing tangles in any such delta diagram of $K$. For links with multiple components, $c_{\Delta}(L)$ is defined similarly using delta-crossing tangles in a link diagram.
	
	\subsection{Other knot invariants}
	We briefly recall definitions of some other standard invariants that will appear. The \emph{crossing number} $c(L)$ of a link $L$ is the minimum number of crossings in any diagram of $L$. The \emph{(orientable) genus} $g(L)$ is the minimum genus of any Seifert surface for $L$. The \emph{triple-crossing number} $c_3(L)$ is the least number of triple intersection points in any diagram of $L$ with only triple points (no ordinary double crossings) \cite{Ada13}. We will also encounter invariants associated with special link diagrams: for example, if $L$ can be represented as the closure of a braid, the \emph{braid index} $b(L)$ is the minimal number of strands in such a braid. In addition, we denote by $P_L(a,z)$ the HOMFLY-PT polynomial of $L$, in variables $a$ and $z$ (with the normalization $P_{\text{unknot}}(a,z)=1$). We write $\deg_z P_L(a,z)$ for the degree of $P_L$ in the $z$-variable.

	\noindent
	\begin{figure}
		\begin{tikzcd}[column sep=1.07cm, row sep=1.5cm]
			2b-2\arrow[<-, bend left=42, pos=.2,swap]{ddrrrrrrrr}{\color{red} 14}&&[-0.5cm]2br-2\arrow[<-, bend left=15, pos=.5]{rrr}{\color{red} 33}\arrow[<-,swap]{ll}{\color{red} 22}&[-0.5cm]2u_c\arrow[<-, bend left=0]{dr}{\color{red} 31}&[-0.5cm]&[-0.5cm]2a\arrow[<-, bend left=10,pos=0.1]{ddrrr}{\color{red} 24}\arrow[->, bend right=30, swap]{dl}{\color{red} 32}&[-0.5cm]&[-0.5cm]&[-0.5cm]\\[-0.5cm]spP_v\arrow[<-, bend right=0, pos=.6]{u}{\color{red} 23}&[-0.5cm]{2cl_4}\arrow[<-]{r}{\color{red} 12}&[-0.5cm]{2u_s}\arrow[<-]{r}{\color{red} 25}&[-0.5cm]2u^*_r\arrow[<-, bend left=2,pos=0.4]{r}{\color{red} 35}&[-0.5cm]2u\arrow[<-, bend left=2]{r}{\color{red} 29}\arrow[<-, bend left=40,swap]{rr}{\color{red} 10}&[-0.5cm]2cl&[-0.5cm]tr\arrow[<-, bend left=0,swap]{drr}{\color{red} 26}\\2g_4\arrow[<-,pos=.5]{r}{\color{red} 19}\arrow[<-, bend left=32,pos=.9]{ur}{\color{red} 11}\arrow[<-]{dr}{\color{red} 21}&[-0.5cm]2g_r \arrow[<-, swap]{r}{\color{red} 18} \arrow[<-, bend left=5,swap]{urr}{\color{red} 30} &[-0.5cm] ul_b \arrow[<-, pos=.4]{r}{\color{red} 17} &[-0.5cm] u_b\arrow[<-, bend left=10,swap]{ur}{\color{red} 15}\arrow[<-, swap]{r}{\color{red} 16}&[-0.5cm]2g \arrow[<-, pos=.7,swap]{dl}{\color{red} 27} \arrow[<-, pos=.7]{r}{\color{red} 6}\arrow[<-, bend right=18]{ur}{\color{red} 4}&[-0.5cm] 2g_f \arrow[<-]{r}{\color{red} 5}&[-0.5cm] 2g_c \arrow[<-]{r}{\color{red} 3}\arrow[<-]{u}[bend left=30,pos=.4,]{\color{red} 28} &[-0.5cm] c_3 \arrow[<-]{r}[bend left=30,swap]{\color{red} 2}\arrow[<-]{d}{\color{red} 36}&[-0.5cm]c\\
			{|\sigma|}\arrow[<-,swap, pos=.4,swap]{u}{\color{red} 13}&[-0.5cm] g_{ds} \arrow[<-, swap]{ur}{\color{red} 20}&[-0.5cm]&[-0.5cm]td\arrow[<-, bend right=35]{rrrr}{\color{red} 43} \arrow[<-, bend right=77,swap]{urrrrr}{\color{red} 1}\arrow[->, bend right=20, pos=.6,swap]{rr}{\color{red} 39}&[-0.5cm]sp\Delta_t\arrow[<-, pos=.6]{r}{\color{red} 40}\arrow[<-, pos=.4]{u}{\color{red} 8}&[-0.5cm]degP_z\arrow[<-, pos=.6,swap]{ur}{\color{red} 7}&[-0.5cm]\lceil\nicefrac{spV_t}{2}\rceil\arrow[<-, bend right=8, swap]{ur}{\color{red} 9}&[-0.5cm]2c_{\Delta}\\
			{2|\tau|}\arrow[<-, pos=.3, bend right=18,swap]{uu}{\color{red} 42}&{|s|}\arrow[<-,swap, pos=.2]{uul}{\color{red} 41}&spF_a\arrow[<-, bend left=35]{rr}{\color{red} 37}&m-1\arrow[<-,swap, pos=.6]{r}{\color{red} 38}&\alpha-2\arrow[<-, bend right=68, pos=.2]{uurrrr}{\color{red} 34}&&&&\\
		\end{tikzcd}
		\caption{Graph of selected inequalities between invariants}
		\label{g1}
	\end{figure}
	

	\section{Inequalities between knot invariants}\label{s3}

	We present the graphs of selected invariant inequalities for any knot $K\hookrightarrow \mathbb{S}^3$. In the diagram, shown in Figure \ref{g1}, an arrow $\rightarrow$ means the relation $\geq$, the invariants are: 
	\ \\
	$\color{blue} c$ -  the crossing number,
	$\color{blue} c_3$ -  triple crossing index,
	$\color{blue} c_{\Delta}$ -  the delta-crossing number,
	$\color{blue} g$ - the (Seifert) three-genus,
	$\color{blue} g_f$ - the free genus,
	$\color{blue} g_c$ - the canonical genus,
	$\color{blue} u$ - the unknotting number,
	$\color{blue} u_b$ - the band-unknotting number,
	$\color{blue} ul_b$ - the band-unlinking number,
	$\color{blue} g_4$ - the slice genus,
	$\color{blue} g_r$ - the ribbon slice genus,
	$\color{blue} g_{ds}$ - the doubly slice genus, 
	$\color{blue} \sigma$ - the knot signature,
	$\color{blue} sp\Delta_t$ - the span of the Alexander polynomial $\Delta(t)$,
	$\color{blue} spV_t$ - the span of the Jones polynomial $V(t)$,
	$\color{blue} degP_z$ - the $z$-degree of the HOMFLYPT polynomial $P(v,z)$,
	$\color{blue} spP_v$ - the $a$-spread of the HOMFLYPT polynomial $P_K(a,z)$,
	$\color{blue} spF_a$ - the $a$-spread of the Kauffman polynomial $F_K(a,z)$,
	$\color{blue} cl_4$ - the 4D clasp number,
	$\color{blue} cl$ - the clasp number,
	$\color{blue} u_s$ - the slicing number,
	$\color{blue} td$ - the skein tree depth,
	$\color{blue} tr$ - the trivializing number,
	$\color{blue} u_c$ - the concordance unknotting number,
	$\color{blue} u^*_r$ - the weak ribbon unknotting number, 
	$\color{blue} br$ - the bridge number,
	$\color{blue} b$ - the braid index,
	$\color{blue} \alpha$ - the arc index,
	$\color{blue} m$ - the mosaic number,
	$\color{blue} a$ - the ascending number,
	$\color{blue} \tau$ - the Ozsvath-Szabo's tau-invariant,
	$\color{blue} s$ - the Rasmussen's s-invariant.

	\par 
	The more precise definitions of the invariants and the inequality sources are as follows.\\ 
	Relation:  
	{\color{red} 1} in \cite{Cro89},
	{\color{red} 2, 9} in \cite{Ada13},
	{\color{red} 3} in \cite{Jab20},
	{\color{red} 5} in \cite{KobKob96},
	{\color{red} 6} in \cite{Mor87},
	{\color{red} 7} in \cite{Mor86},
	{\color{red} 8} in \cite{Gil82},
	{\color{red} 4, 11, 12, 19, 25, 29, 30, 35} in \cite{Shi74},
	{\color{red} 26} in \cite{Han14},
	{\color{red} 13} in \cite{KauTay76},
	{\color{red} 16, 17, 18} in \cite{JMZ20},
	{\color{red} 15} in \cite{HNT90},
	{\color{red} 20} in \cite{Mcd19},
	{\color{red} 21} in \cite{LivMei15},
	{\color{red} 31} in \cite{OweStr16},
	{\color{red} 27} in \cite{SchTho89},
	{\color{red} 10, 28} in \cite{Hetal11},
	{\color{red} 36} in \cite{Jab23},
	{\color{red} 38} in \cite{LHLO14},
	{\color{red} 37} in \cite{MorBel98},
	{\color{red} 34} in \cite{BaePark00},
	{\color{red} 24, 32, 33} in \cite{Oza10},
	{\color{red} 23} in \cite{FW85, Mor86},
	{\color{red} 22} in \cite{Yam87},
	{\color{red} 14} in \cite{Ohy93},
	{\color{red} 39} in this article.
	{\color{red} 40} is immediate,
	{\color{red} 41} in \cite{Ras10},
	{\color{red} 42} in \cite{OzsSza03},
	{\color{red} 43} in this article.

	\section{Inequality relating $\deg_z P$ and $td$}\label{s4}
	We first summarize several known bounds relating the skein tree depth to other invariants of $L$. The following gives a general lower and upper bound for $td(L)$ in terms of topological invariants of $L$:
	
	\begin{theorem}[{\cite{Cro89,SchTho89}}]\label{thm:genus-crossing-bounds}
		For any nontrivial link $L$, 
		\[ 2\,g(L) + r(L) - 1 \ \le\ td(L)\ \le\ c(L) - 1, \] 
		where $g(L)$ is the oriented genus of $L$, $r(L)$ is the number of components of $L$, and $c(L)$ is the crossing number of $L$.
	\end{theorem}
	
	Another known result gives the exact skein tree depth for links that are closures of certain pure braids:
	
	\begin{theorem}[{\cite{KKO15}}]\label{thm:positive-braid}
		Suppose $L$ is a (non-split) link that can be represented as the closure of a braid $\beta$ on $p$ strands, where $\beta$ is a word of length $k$ consisting entirely of either all positive or all negative braid generators (so $\beta$ is a positive or negative braid). Then 
		\[ td(L) = k - p + 1.\] 
		In particular, $td(L)$ is determined exactly by the braid word length in this case.
	\end{theorem}
	
	We also have a general lower bound for $td(L)$ coming from the HOMFLY-PT polynomial. Specifically, the skein relation for HOMFLY-PT implies that each skein move can increase the $z$-degree of the polynomial by at most 1. Since an $r$-component unlink has $\deg_z P_{\text{unlink}} \le 0$, it follows that enough skein moves must be performed to climb from degree 0 up to the $z$-degree of $P_L$. This intuitive argument is made precise in the following proposition:
	
	\begin{proposition}\label{prop:degPz-bound}
		For any nontrivial link $L$, 
		\[ td(L)\ \ge\ \deg_z P_L(a,z). \]
	\end{proposition}
	
	\begin{proof}
		The HOMFLY-PT skein relation is given by 
		\[a^{-1} P_{L_+}(a,z)\ -\ a\,P_{L_-}(a,z)\ =\ z\,P_{L_0}(a,z).\] 
		For an $r$-component unlink $U_r$, one has $P_{U_r}(a,z) = ((a^{-1}-a)/z)^{\,r-1}$, which is a polynomial of $z$-degree $0$ (indeed a monomial $z^{-(r-1)}$). Each single skein resolution can increase the $z$-degree of the polynomial by at most $1$, because of the factor of $z$ multiplying $P_{L_0}(a,z)$ in the skein relation. Therefore, in any skein tree for $L$, the distance from the root $L$ (with polynomial $P_L$ of some $z$-degree $d$) to a leaf $U_r$ (with polynomial of $z$-degree $\le 0$) must be at least $d$. Hence, the minimal depth $td(L)$ cannot be less than $\deg_z P_L(a,z)$.
	\end{proof}
	
	The bound from Proposition~\ref{prop:degPz-bound} can be strictly stronger than the genus bound in Theorem~\ref{thm:genus-crossing-bounds} for certain links. For example, for the knot $K = 11{n42}$ one finds $2\,g(K) = 4$ whereas $\deg_z P_K = 6$, so Proposition~\ref{prop:degPz-bound} gives $td(K)\ge 6$ compared to the weaker $td(K)\ge 4$ from genus. On the other hand, the $\deg_z P$ bound is not universally stronger: there exist knots (e.g. $15{n14891}$, see \cite{Jab23}) for which $2\,g(K)$ exceeds $\deg_z P_K$. In general, one should use the maximum of these lower bounds for a given link. 
	
	We note also that all currently known data are consistent with the conjecture that $td(L)$ is bounded below by the \emph{triple-crossing number} $c_3(L)$ (see Conjecture~\ref{conj:triple} in Section\;\ref{s8}).
	
	\section{Inequality relating $c_{\Delta}$ and $td$}\label{s5}
	We now present an inequality that directly connects the delta-crossing number of a knot to its skein tree depth. This result is new and does not follow from any previous relations among invariants (for instance, it is independent of the inequalities summarized in Section\;\ref{s3}).
	
	\begin{theorem}\label{thm:delta-inequality}
		For any knot $K$, we have 
		\[ 2\,c_{\Delta}(K)\ \ge\ td(K).\]
	\end{theorem}
	
	\begin{proof}
		Let $K$ be a knot, and take a knot diagram $D_1$ realizing $c_{\Delta}(K)=n$, i.e. $D_1$ is a minimal delta-crossing diagram with $n$ delta-crossing tangles. Each delta-crossing tangle in $D_1$ can be one of four types (often labeled $S$, $T$, $U$, or $W$) as classified in \cite{NSS15}, see Figure \ref{STUWb}.

			\begin{figure}[h!t]
			\centering		
			\begin{lpic}[]{./STUWb(11.8cm)}
				
			\end{lpic}
			\caption{Four types of delta-crossing.}
			\label{STUWb}
		\end{figure}

		\begin{figure}[h!t]
			\centering		
			\begin{lpic}[]{./delta_tee(14.5cm)}
				
			\end{lpic}
			\caption{Resolving locally the delta-crossing.}
			\label{delta_tee}
		\end{figure}

		We will construct a skein tree for $K$ by resolving each of these delta-crossings in at most two steps, thereby showing $td(K)\le 2n$.
		
		Without loss of generality, assume the first delta-crossing encountered as we traverse $D_1$ is of type $S$ (the other cases $T$, $U$, $W$ are analogous). A delta-crossing contains six crossing arcs (labeled $a,b,c,d,e,f$); if we perform a certain crossing change within the delta tangle, we can make the diagram locally descending along one arc of the tangle.  It is well-known that if we made the diagram descending by the classical crossing changes then the link diagram would be a diagram of a trivial link. Consider this delta-crossing as one shown in Figure \ref{delta_tee} top (root position). Without loss of generality, we assume that we travel in the $a\to d$ direction. Because it is a delta-crossing, the $a\to d$ arc alternates, we are able to make at most one crossing change on  $a\to d$ to obtain our descending diagram procedure. Let it be the crossing marked by a dashed red circle.
		\par 
		We proceed to produce the two child diagrams in the skein tree as shown. If in the left child, the descending procedure forces us to make the crossing changes in the delta-crossing (outside $a\to d$ arc) then we create from this diagram two more child diagrams as shown in the figure (resolving crossing is marked by a dashed red circle). It is not mandatory, as indicated by the question mark in the figure. In the case of the right child of the root diagram, we always make its skein resolution.
		\par
		We notice now that the first of the bottom four pictures in Figure \ref{delta_tee} is proper to guarantee the resulting final diagram $D_2$ after visiting every delta-crossing to be a diagram of the trivial knot. In the case of the other tree bottom pictures, it can be easily seen that performing Reidemeister moves of type I or II completely removes the crossings inside the delta-crossing producing delta-diagrams with $n-1$ delta-crossings, therefore we can treat each resulting delta-diagram as new $D_1*$ diagram being the root of our procedure, inductively reduces to diagram without any classical crossings, i.e. a diagram of a trivial link.   
		\par
		A diagram after the zeroth resolution of its crossing may change its number of link components, but the orientations on any remaining arcs do not change, so the resulting link has all its remaining delta-crossings (beside the one where the resolution takes place) is still one of four possible delta-crossings $S, T, U, W$ like in the starting (root) knot.
		\par
		At the end, we obtain a skein tree of two levels for each delta-crossing as a local root, therefore $td(K)\leq \text{our tree depth} =2n=2c_{\Delta}(K)$.
	\end{proof}
	
	The proof of Theorem~\ref{thm:delta-inequality} is constructive: it shows how to obtain a skein tree of depth $2\,c_{\Delta}(K)$ by a specific strategy of resolving delta-crossings. In particular, each delta-crossing tangle is handled locally by at most two skein moves. The bound is sharp for many knots. For example, for the knot $12a{58}$ one finds $\deg_z P{12a58} = 8$ \cite{LivMoo23} and $c_{\Delta}(12a58) = 4$ \cite{Jab23}, so our inequalities give $8 \le td(12a58) \le 8$. Hence $td(12a58)=8$. The same reasoning shows $td(K)=8$ for a number of other 12-crossing knots, including $12a99, 12a119, 12a268, 12a281, 12a323, 12a426, 12a435, 12a499, 12a561, 12a629, 12a868,\\ 12a1188, 12n188, 12n326, 12n327, 12n328, 12n416, 12n417, 12n425, 12n426, 12n518, 12n538,\\ 12n591, 12n592, 12n609, 12n703, 12n706.$ These examples were previously unknown and demonstrate instances where the bound in Theorem~\ref{thm:delta-inequality} is exact.
	
	\section{New upper bound from braid representations}\label{s6}
	We now derive a new upper bound for skein tree depth in terms of braid parameters. Let $L$ be an oriented link that can be presented as the closure of a (non-split) braid $\beta$ on $s_\beta$ strands. Let $c_\beta(L)$ denote the number of crossings in the braid word (i.e., the length of the word), and let $c^+_\beta(L)$ and $c^-_\beta(L)$ be the number of positive and negative crossings in that braid word, respectively. Clearly, $c_\beta(L) = c^+_\beta(L) + c^-_\beta(L)$ is at least the crossing number $c(L)$ of $L$, but it could be larger for a non-minimal braid representation.
	
	We have the following inequality for $td(L)$:
	
	\begin{theorem}\label{thm:braid-bound}
		For any nontrivial link $L$,
		\[ td(L)\ \le\ \min_{\beta}\Big\{\,c_{\beta}(L) - s_{\beta}(L) + 1 \,+\, \min\big(c^+_{\beta}(L),\;c^-_{\beta}(L)\big)\Big\}, \]
		where the minimum is taken over all braid representations $\beta$ of $L$. 
	\end{theorem}
	
	In simpler terms, for any fixed braid presentation of $L$, the quantity $c_{\beta}(L) - s_{\beta}(L) + 1 + \min(c^+_{\beta},c^-_{\beta})$ is an upper bound on $td(L)$, and one gets the sharpest bound by minimizing this over all braids for $L$.
	
	\begin{proof}
It is sufficient to show an inequality for the fixed, arbitrary braid $\beta$ such that its braid closure is $L$. Without loss of generality assume that $c_{\beta}^+(L) \geq c_{\beta}^-(L)$. We now place the braid closure of $\beta$ at the root of the skein tree. Now, we resolve all possible negative crossings by skein relation, such that the zeroth resolution is still a non-split braid. We do this step by step, fixing always the link $L$ neighborhoods of each crossing. In each such step, the negative and the zeroth skein resolutions of an oriented diagram have one less negative crossing and either a fixed number of positive crossings or the same number of positive crossings.
\par
We are left now, with the tree leaves on at most level equal to $c_{\beta}^-(L)$ from the root, with the possible negative crossings such that making the zeroth resolution on each of them produces a split link. In this case, they are so-called nugatory crossings, and we just perform a full twist on part of a link diagram from one side of each of the nugatory crossings, switching them into a positive crossing.
\par
Now all diagrams are diagrams of positive braids $\beta_{k}$ corresponding to links $L_{k}$ such that each of them satisfies inequality $c_{\beta_{k}}(L_{k})-s_{\beta_{k}}(L_{k})+1 \leq c_{\beta}(L)-s_{\beta}(L)+1$ because $s_{\beta}(L)=s_{\beta_{k}}(L_{k})$ as they may be split links. Applying now Theorem\;\ref{thm:positive-braid} to $\beta_{k}$ we obtain the skein tree depth of each closure of $\beta_{k}$ that is at most $c_{\beta_{k}}(L_{k})-s_{\beta_{k}}(L_{k})+1$, hence our desired inequality in the theorem is true by adding the path of length equal at most $c_{\beta}^-(L)$ to the root $\beta$, making the whole tree of depth at most $c_{\beta}(L)-s_{\beta}(L)+1+c_{\beta}^-(L)$.
	\end{proof}
	
	The bound in Theorem~\ref{thm:braid-bound} can outperform the general $td(L)\le c(L)-1$ bound from Theorem~\ref{thm:genus-crossing-bounds} for some knots and links (that are moreover neither positive nor negative braids to apply Theorem\;\ref{thm:positive-braid}). Consider for example the knot $K11n183$, and the closure of one of its corresponding braids $\sigma_1^{-1}\sigma_2\sigma_1^{-1}\sigma_3^{-1}\sigma_2^{-2} \sigma_1^{-1}\sigma_3^{-1}\sigma_2^{-2}\sigma_3^{-1}$. This gives us a new upper bound on the skein tree depth equal to at most $9$. The other such knot, with a stronger upper bound on the skein tree depth, is $K12n163$ and for the non-knot example, we have $L11n443\{0,1,1\}$. The latter link has the Conway polynomial equal $-z^3+4z^5+z^7$ (see \cite{LivMoo23}). We see that not all coefficients are either positive or negative, so it is neither positive nor negative braid by \cite{Cro89}.
	
	\section{Examples and tabulations}\label{s7}
	In this section, we describe the computational determination of skein tree depth for many knots and links. Tables~\ref{table1}--\ref{table3} list the values (or, in some cases, narrow ranges) of $td(L)$ for all prime knots and links up to certain complexity. Table~\ref{table1} covers all oriented prime knots and links with crossing number $c(L)\le 8$. Table~\ref{table2} lists additional knots and links with $c(L)=9$ for which we were able to compute an exact $td$ value. Table~\ref{table3} lists selected prime knots with $c(L)\ge 10$ for which we determined $td$ exactly. Together, these tables significantly extend the range of links with known skein tree depth.
	
	The computations were performed by a recursive algorithm implementing the skein resolution process with pruning optimizations. We considered each prime knot/link up to 8 crossings (using the standard Dowker-Thistlethwaite link listings up to 8 crossings) and many with higher crossings. We treated links up to mirror image, since $td(L)$ is the same for a link and its mirror. For links with a specified orientation convention (involving braces in the notation), we followed the orientation conventions given in \cite{LivMoo23}.
	
	The algorithm uses several strategies to manage the potentially huge skein trees:
	\begin{itemize}\setlength\itemsep{0pt}
		\item We start from a basis of all reduced link diagrams up to 8 crossings provided in \cite{Jab19} and only consider those as initial diagrams, to avoid redundant exploration of many equivalent diagrams.
		\item We utilize known values of $td(L)$ for certain links to terminate branches early. For instance, all links with $td(L)\le 2$, and values for those links are taken as known (these known cases are indicated in bold in Table~\ref{table1}). Similarly, links that fall under Theorem~\ref{thm:positive-braid} (closures of pure $\pm$ braids) have $td$ computed by the formula and can be used to prune search branches.
		\item We apply the inequalities from Theorem~\ref{thm:genus-crossing-bounds} (and Proposition~\ref{prop:degPz-bound}) as pruning conditions. If at any point a partial resolution yields a diagram $D$ for which the known lower bound on $td$ (using genus or polynomial degree) is greater than the remaining depth budget in our search, we can stop exploring that branch. In some cases, for multi-component links, we used the breadth (span) of the one-variable Alexander polynomial in place of $2g + r - 1$ as a lower bound, since for split links the genus bound can be less effective.
	\end{itemize}
	
	Despite these optimizations, certain links present a challenge: in a few instances, the recursive algorithm falls into a potentially infinite loop of generating equivalent diagrams under different sequences of resolutions. When such a situation was detected (the algorithm revisiting a diagram it had seen before in the same branch), we used the lower bound from Theorem~\ref{thm:genus-crossing-bounds} (or an Alexander polynomial bound for links) to break the loop, assigning the unresolved branch a provisional depth equal to that lower bound. This ensured termination of the search.
	
	The compiled results are displayed in Table~\ref{table1} (for $c\le 8$), Table~\ref{table2} (selected $c=9$ cases), and Table~\ref{table3} (selected $c\ge 10$ cases). In these tables, a single integer denotes the exact value of $td(L)$; a range $[m,n]$ denotes that $td(L)$ lies between $m$ and $n$ (inclusive), and we could not narrow it down further with our current computations. Boldface entries in Table~\ref{table1} indicate links for which $td$ was known from prior results (specifically, those with $td\le 2$ from \cite{SchTho89} and those covered by Theorem~\ref{thm:positive-braid}).
	
	One noteworthy pattern is that all computed values support the conjecture that $td(L) \ge c_3(L)$ for any link $L$ (see Conjecture~\ref{conj:triple} in the next section). In fact, in our tables, $td(L)$ often equals $c_3(L)$ for the entries where $c_3$ is known \cite{JabTro20, Jab23a}. Additionally, we observe that our new delta-crossing inequality (Theorem~\ref{thm:delta-inequality}) frequently provides the exact value of $td$ when used in combination with Proposition~\ref{prop:degPz-bound}, as illustrated by the examples of 12-crossing knots discussed earlier.
	
\begin{footnotesize}
	\renewcommand{\arraystretch}{1.25}
	\begin{center}
		
		\begin{longtable}[ht]{r||l}
			\caption{Knots and links and the skein tree depth\label{table1}}\\
			$td(L)$	&  names of prime knots or links $L$ with $c(L)\leq 8$\\
			\hline
			\endfirsthead
			\multicolumn{2}{c}
			{\tablename\ \thetable\ -- \textit{Continued from previous page}} \\
			$td(L)$	 & names of prime knots or links $L$ with $c(L)\leq 8$\\
			\hline
			\endhead
			\hline \multicolumn{2}{r}{\textit{Continued on next page}} \\
			\endfoot
			\hline
			\endlastfoot
			$1$&$\bf L2a1$.\\
			\hline
			$2$&$\bf K3a1$, $\bf K4a1$, $\bf L4a1\{0\}$.\\
			\hline
			$3$&$K5a1$, $K6a3$, $\bf L4a1\{1\}$, $L5a1$, $L6a1\{0\}$, $L6a3\{1\}$, $L6a5\{0,0\}$, $L6n1\{0,0\}$, $L6n1\{1,0\}$,
			$L6n1\{1,1\}$.
			\\
			\hline
			$[3, 4]$&$K7a4$, $K7a6$, $K8a11$, $K8a18$, $L6a1\{1\}$, $L6a2$, $L7a2\{0\}$, $L7a4$, $L7a5\{0\}$, $L7a6\{1\}$, $L7n1\{1\}$, $L7n2$,
			\\
			& 
			$L8a3\{0\}$, $L8a6\{0\}$, $L8a11\{1\}$, $L8a14\{1\}$, $L8a18\{0,1\}$, $L8a21\{0,0,0\}$, $L8n1\{1\}$, $L8n2$, $L8n3\{1,0\}$, 
			\\
			&
			$L8n3\{0,1\}$, $L8n4\{0,1\}$, $L8n7\{0,0,0\}$, $L8n7\{1,0,0\}$, $L8n7\{1,0,1\}$, $L8n7\{1,1,1\}$, $L8n8\{0,0,0\}$, 
			\\
			&
			$L8n8\{1,0,0\}$, $L8n8\{0,1,0\}$, $L8n8\{1,1,0\}$, $L8n8\{0,0,1\}$, $L8n8\{1,1,1\}$.
			\\
			\hline
			$4$&$\bf K5a2$, $K6a1$, $K6a2$, $K7a1$, $K7a2$, $K8a5$, $L6a4$, $L6a5\{0,1\}$, 
			
			$L6a5\{1,0\}$, $L6a5\{1,1\}$, $\bf  L6n1\{0,1\}$, 
			\\
			&
			$L7a7\{0,0\}$, $L7a7\{0,1\}$, $L7a7\{1,0\}$, $L8a15\{0,0\}$, $L8a20\{0,0\}$, $L8n5\{1,0\}$, $L8n5\{0,1\}$, $L8n6\{0,0\}$.
			\\
			\hline
			$[4,5]$&$K7a3$, $K7a5$, $K8a1$, $K8a2$, $K8a4$, $K8a7$, $K8a9$, $K8a10$, $K8a17$,$K8n1$, $K8n2$, $L7a7\{1,1\}$,  
			\\
			&
			$L8a15\{1,0\}$, $L8a15\{0,1\}$, $L8a15\{1,1\}$, $L8a16\{1,0\}$, $L8a16\{0,1\}$, $L8a17\{0,0\}$, $L8a17\{1,0\}$, 
			\\
			&
			$L8a17\{0,1\}$, $L8a18\{1,0\}$, $L8a19\{1,0\}$, $L8a19\{0,1\}$, $L8a20\{0,1\}$, $L8a20\{1,1\}$,
			$L8n3\{1,1\}$, 	
			\\
			&
			$L8n4\{0,0\}$, 	$L8n4\{1,0\}$, 	$L8n4\{1,1\}$, $L8n5\{0,0\}$, 	$L8n5\{1,1\}$, $L8n6\{0,1\}$, 	$L8n6\{1,1\}$.
			\\
			\hline
			$5$&$\bf L6a3\{0\}$, $L7a1$, $L7a2\{1\}$, $L7a3$, $L7a5\{1\}$, $L7a6\{0\}$, $\bf L7n1\{0\}$, $L8a1$, $L8a2$, $L8a3\{1\}$,  $L8a4$, $L8a5\{0\}$,
			\\
			&
			
			$L8a7\{1\}$, $L8a8$, $L8a9$, $L8a10\{1\}$,
			$L8a21\{1,0,0\}$, $L8a21\{0,1,0\}$, $L8a21\{0,0,1\}$, $L8a21\{1,0,1\}$,  
			\\
			&
			$L8a21\{0,1,1\}$, $L8a21\{1,1,1\}$, $L8n1\{0\}$, $L8n7\{0,1,0\}$, $L8n7\{1,1,0\}$, $L8n7\{0,0,1\}$, $L8n7\{0,1,1\}$, 
			\\
			&
			$\bf  L8n8\{1,0,1\}$, $\bf L8n8\{0,1,1\}$.
			\\
			\hline
			$[3, 4, 5]$& $L8a6\{1\}$,  $L8a7\{0\}$, $L8a10\{0\}$, $L8a12\{1\}$, $L8a13\{0\}$.
			\\
			\hline
			$[5, 6]$& $L8a5\{1\}$, $L8a11\{0\}$, $L8a12\{0\}$, $L8a13\{1\}$, $L8a21\{1,1,0\}$.
			\\
			\hline
			$6$&$\bf K7a7$, $K8a3$, $K8a6$, $K8a8$,  $K8a12$, $K8a13$ ,$K8a14$, $K8a15$, $K8a16$, $\bf K8n3$, $L8a16\{0,0\}$, $L8a16\{1,1\}$,
			\\
			&
			$L8a17\{1,1\}$, $L8a18\{0,0\}$, $L8a18\{1,1\}$, $L8a19\{0,0\}$, $L8a19\{1,1\}$, $L8a20\{1,0\}$, $\bf L8n3\{0,0\}$,  
			$\bf L8n6\{1,0\}$.
			\\
			\hline
			$7$&$\bf L8a14\{0\}$.
			\\
			\hline
		\end{longtable}
	\end{center}
	
\end{footnotesize}

	
	\begin{footnotesize}
		\renewcommand{\arraystretch}{1.25}
		\begin{center}
			
			\begin{longtable}[ht]{r||l}
				\caption{Knots and links and the skein tree depth\label{table2}}\\
				$td(L)$	&  names of prime knots or links $L$ with $c(L) = 9$\\
				\hline
				\endfirsthead
				\multicolumn{2}{c}
				{\tablename\ \thetable\ -- \textit{Continued from previous page}} \\
				$td(L)$	 & names of prime knots or links $L$ with $c(L) = 9$\\
				\hline
				\endhead
				\hline \multicolumn{2}{r}{\textit{Continued on next page}} \\
				\endfoot
				\hline
				\endlastfoot
				$4$& $K9n5$, $K9n6$.
				
				\\
				\hline
				$5$& $L9a5\{0\}$, $L9a8\{0\}$, $L9a8\{1\}$, $L9a11\{0\}$, $L9a16\{0\}$, $L9a26\{1\}$, $L9a27\{0\}$, $L9a33\{0\}$, $L9a42\{1\}$,
				\\
				&
				$L9a55\{0,0,0\}$, $L9a55\{0,1,0\}$, $L9a55\{0,0,1\}$, $L9a55\{1,0,1\}$, $L9n8\{0\}$, $L9n8\{1\}$, $L9n9\{0\}$, $L9n10\{0\}$,  
				\\
				&
				$L9n11\{0\}$, $L9n19\{0\}$, $L9n19\{1\}$.
				
				\\
				\hline
				$6$& $K9a1$, $K9a2$, $K9a5$, $K9a6$, $K9a7$, $K9a9$, $K9a11$, $K9a12$, $K9a13$, $K9a14$, $K9a15$, $K9a19$, $K9a20$, 
				\\
				& 
				$K9a28$, $K9a31$, $K9a37$, $K9n3$, $K9n7$, $L9a43\{1,0\}$, $L9a43\{0,1\}$, $L9a43\{1,1\}$, $L9a44\{0,0\}$,
				\\
				&
				$L9a44\{1,0\}$, $L9a44\{0,1\}$, $L9a46\{0,0\}$, $L9a46\{1,0\}$, $L9a46\{0,1\}$, $L9a46\{1,1\}$, $L9a47\{1,0\}$,
				\\
				&
				$L9a47\{0,1\}$, $L9a47\{1,1\}$, $L9a48\{0,0\}$, $L9a49\{1,0\}$, $L9a49\{0,1\}$, $L9a50\{0,0\}$, $L9a50\{1,0\}$,
				\\
				&
				$L9a50\{1,1\}$, $L9a51\{0,0\}$, $L9a51\{1,0\}$, $L9a51\{1,1\}$, $L9a52\{1,0\}$, $L9a52\{1,1\}$,  $L9a53\{0,0\}$,
				\\
				&
				$L9a53\{1,0\}$, $L9a53\{0,1\}$, $L9a53\{1,1\}$, $L9a54\{0,0\}$, $L9a54\{1,0\}$, $L9a54\{0,1\}$, $L9a54\{1,1\}$,
				\\
				&
				$L9n20\{0,1\}$, $L9n21\{0,1\}$, $L9n22\{1,0\}$, $L9n22\{0,1\}$, $L9n22\{1,1\}$, $L9n23\{0,0\}$, $L9n24\{1,0\}$,
				\\
				&
				$L9n26\{1,0\}$, $L9n26\{1,1\}$, $L9n28\{0,0\}$, $L9n28\{1,1\}$.
				\\
				\hline
				$7$& $L9a2\{0\}$, $L9a2\{1\}$, $L9a6\{1\}$, $L9a9\{0\}$, $L9a9\{1\}$, $L9a12\{1\}$, $L9a14\{0\}$, $L9a14\{1\}$, $L9a20\{0\}$, $L9a21\{0\}$,
				\\
				& 
				$L9a22\{0\}$, $L9a24\{1\}$, $L9a28\{0\}$, $L9a29\{0\}$, $L9a31\{0\}$, $L9a32\{1\}$, $L9a36\{0\}$, $L9a38\{0\}$, $L9a39\{0\}$,
				\\
				&
				$L9a41\{0\}$, $L9a42\{0\}$, $L9n4\{0\}$, $L9n12\{1\}$, $L9n15\{0\}$, $L9n18\{0\}$.
				\\
				\hline
				$8$&$K9a41$.
				\\
				\hline
			\end{longtable}
		\end{center}
		
	\end{footnotesize}

	\begin{footnotesize}
		\renewcommand{\arraystretch}{1.25}
		\begin{center}
			
			\begin{longtable}[ht]{r||l}
				\caption{Knots and the skein tree depth\label{table3}}\\
				$td(K)$	&  names of prime knots $K$ with $c(K)\geq 10$\\
				\hline
				\endfirsthead
				\multicolumn{2}{c}
				{\tablename\ \thetable\ -- \textit{Continued from previous page}} \\
				$td(K)$	 & names of prime knots $K$ with $c(K)\geq 10$\\
				\hline
				\endhead
				\hline \multicolumn{2}{r}{\textit{Continued on next page}} \\
				\endfoot
				\hline
				\endlastfoot
				$6$& $10a1$,
				$10a2$,
				$10a3$,
				$10a7$,
				$10a10$,
				$10a11$,
				$10a17$,
				$10a21$,
				$10a22$,
				$10a24$,
				$10a25$,
				$10a27$,\\&
				$10a31$,
				$10a32$,
				$10a35$,
				$10a38$,
				$10a52$,
				$10a53$,
				$10a66$,
				$10a72$,
				$10a94$,
				$10a100$,
				$10n1$,\\&
				$10n7$,
				$10n8$,
				$10n9$,
				$10n31$,
				$10n32$,
				$10n33$,
				$10n35$,
				$10n41$, 
				$11a5$,
				$11a17$,
				$11a42$,
				$11a51$,\\&
				$11a96$,
				$11a121$,
				$11a128$,
				$11a159$,
				$11a209$,
				$11a218$,
				$11a228$,
				$11n7$,
				$11n8$,
				$11n31$,
				$11n32$,\\&
				$11n33$,
				$11n66$,
				$11n98$,
				$11n103$,
				$11n115$,
				$11n121$,
				$11n124$,
				$11n131$,
				$11n135$,
				$11n137$,\\&
				$11n163$,
				$11n164$,
				$11n165$,
				$11n168$,
				$11n179$,
				$11n183$

				\\
				\hline
				$8$& $10a15$,
				$10a41$,
				$10a56$,
				$10a59$,
				$10a76$,
				$10a78$,
				$10a79$,
				$10a81$,
				$10a83$,
				$10a86$,
				$10a88$,\\&
				$10a91$,
				$10a93$,
				$10a95$,
				$10a103$,
				$10a104$,
				$10a106$,
				$10a107$,
				$10a110$,
				$10a118$,
				$10a120$,\\&
				$10a121$,
				$10a122$,
				$10n21$,
				$10n27$,
				$10n36$, 
				$11a3$,
				$11a7$,
				$11a9$,
				$11a14$,
				$11a15$,
				$11a19$,
				$11a22$,\\&
				$11a24$,
				$11a25$,
				$11a26$,
				$11a28$,
				$11a33$,
				$11a34$,
				$11a35$,
				$11a40$,
				$11a44$,
				$11a47$,
				$11a53$,\\&
				$11a55$,
				$11a57$,
				$11a62$,
				$11a66$,
				$11a68$,
				$11a71$,
				$11a72$,
				$11a73$,
				$11a74$,
				$11a76$,
				$11a79$,\\&
				$11a80$,
				$11a81$,
				$11a82$,
				$11a83$,
				$11a86$,
				$11a88$,
				$11a92$,
				$11a99$,
				$11a106$,
				$11a108$,
				$11a109$,\\&
				$11a112$,
				$11a113$,
				$11a125$,
				$11a126$,
				$11a127$,
				$11a129$,
				$11a131$,
				$11a139$,
				$11a142$,
				$11a146$,\\&
				$11a147$,
				$11a151$,
				$11a156$,
				$11a157$,
				$11a158$,
				$11a160$,
				$11a162$,
				$11a163$,
				$11a164$,
				$11a170$,\\&
				$11a171$,
				$11a174$,
				$11a175$,
				$11a176$,
				$11a177$,
				$11a179$,
				$11a180$,
				$11a182$,
				$11a184$,
				$11a189$,\\&
				$11a194$,
				$11a196$,
				$11a203$,
				$11a206$,
				$11a215$,
				$11a216$,
				$11a217$,
				$11a221$,
				$11a223$,
				$11a231$,\\&
				$11a232$,
				$11a233$,
				$11a239$,
				$11a248$,
				$11a251$,
				$11a252$,
				$11a253$,
				$11a254$,
				$11a255$,
				$11a257$,\\&
				$11a259$,
				$11a261$,
				$11a264$,
				$11a266$,
				$11a267$,
				$11a268$,
				$11a269$,
				$11a274$,
				$11a277$,
				$11a281$,\\&
				$11a282$,
				$11a284$,
				$11a286$,
				$11a287$,
				$11a288$,
				$11a289$,
				$11a293$,
				$11a300$,
				$11a301$,
				$11a302$,\\&
				$11a305$,
				$11a306$,
				$11a308$,
				$11a314$,
				$11a315$,
				$11a316$,
				$11a326$,
				$11a330$,
				$11a332$,
				$11a346$,\\&
				$11a348$,
				$11a350$,
				$11a351$,
				$11n9$,
				$11n13$,
				$11n23$,
				$11n27$,
				$11n36$,
				$11n41$,
				$11n44$,
				$11n47$,\\&
				$11n57$,
				$11n60$,
				$11n61$,
				$11n76$,
				$11n77$,
				$11n78$,
				$11n81$,
				$11n88$,
				$11n104$,
				$11n107$,
				$11n120$,\\&
				$11n133$,
				$11n147$,
				$11n148$,
				$11n149$,
				$11n153$,
				$11n158$,
				$11n166$,
				$11n173$,
				$11n177$,
				$11n182$

				\\
				\hline
				$10$& $11a367$, $12n242$, $12n472$, $12n574$, $12n679$, $12n688$, $12n725$, $12n888$

			\end{longtable}
		\end{center}
		
	\end{footnotesize}
	
	\section{Open problems and conjectures}\label{s8}
	We conclude with some open questions arising from this work, concerning the relationships between skein tree depth and other invariants:
	
	\begin{question}\label{ques:minimal-diagram}
		Is the skein tree depth $td(L)$ of an oriented knot or link $L$ always achieved on a minimal crossing diagram of $L$? 
	\end{question}
	
	\begin{conjecture}\label{conj:triple}
		For any link $L$, we have $td(L) \ge c_3(L)$.
	\end{conjecture}
	
	Conjecture~\ref{conj:triple} was stated (as Conjecture~4.2) in \cite{Jab23} and is supported by all currently known data.
	
	Finally, we present a list of chosen relations, involving the delta-crossing number, that remain conjectural. It does not follow immediately (by transitivity) from relations in Figure \ref{g1}, and there are no data in \cite{LivMoo23} that contradict them.

	\begin{question}
		Do the following relations hold for any knot $K$? 
		
		\begin{enumerate}
			\item $\;\;cl(K) \leq c_{\Delta}(K)$,
			\item $\;\;tr(K) \leq 2c_{\Delta}(K)$,
			\item $\;\;a(K) \leq c_{\Delta}(K)$.
			
		\end{enumerate}
	\end{question}

\end{document}